\RequirePackage{fix-cm}
\documentclass[twocolumn]{svjour3}          
\smartqed  
\usepackage{graphicx}
\usepackage{color}
\usepackage{amsmath}
\usepackage{amsfonts}
\usepackage{float}
\usepackage{comment}
\usepackage{algorithm} 
\usepackage{algorithmic}  
\usepackage[algo2e]{algorithm2e} 
\usepackage{lipsum}
\usepackage{tabularx}

\journalname{Health Care Management Science}

\usepackage{geometry}
\begin{document}
\sloppy

\title{Managing access to primary care clinics using robust scheduling templates
}	
\titlerunning{Managing access to primary care clinics}        

\author{Sina Faridimehr \and
	Saravanan Venkatachalam \and 
	Ratna Babu Chinnam
}


\institute{Sina Faridimehr \at
	{WarnerMedia}\\
	1050 Techwood Dr NW, Atlanta, GA 30318\\
	Tel.: +1(404)-827-1500\\
	\email{sina.faridimehr@gmail.com}
	\and
	Saravanan Venkatachalam \at
	Department of Industrial and Systems Engineering\\
	Wayne State University\\
	Detroit, MI 48201\\
	Tel.: +1(313)-577-1821\\
	\email{saravanan.v@wayne.edu}
	\and
	Ratna Babu Chinnam \at
	Department of Industrial and Systems Engineering\\
	Wayne State University\\
	Detroit, MI 48201\\
	Tel.: +1(313)-577-4846\\
	\email{Ratna.Chinnam@wayne.edu}
}

\date{Received: date / Accepted: date}

\maketitle

\begin{abstract}
	An important challenge confronting healthcare is the effective management of \textit{access} to primary care. Robust appointment scheduling policies/templates can help strike an effective balance between the lead-time to an appointment (a.k.a. \textit{indirect} waiting time, measuring the difference between a patient's desired and actual appointment dates) and waiting times at the clinic on the day of the appointment (a.k.a. \textit{direct} waiting time). We propose methods for identifying effective appointment scheduling templates using a two-stage stochastic mixed-integer linear program model. The model embeds simulation for accurate evaluation of direct waiting times and uses sample average approximation method for computational efficiency. The model accounts for patients' no-show behaviors, provider availability, overbooking, demand uncertainty, and overtime constraints. The model allows the scheduling templates to be potentially updated at regular intervals while minimizing the patient expected waiting times and balancing provider utilization. Proposed methods are validated using data from the U.S. Department of Veterans Affairs (VA) primary care clinics. 
	\keywords{appointment scheduling \and patient no-show \and direct waiting \and indirect waiting \and stochastic programming \and two-stage model \and sample average approximation \and simulation }
\end{abstract}

\section{Introduction}
\label{intro}
\begin{figure}[!t]
	\centering
	\includegraphics[scale = 0.5]{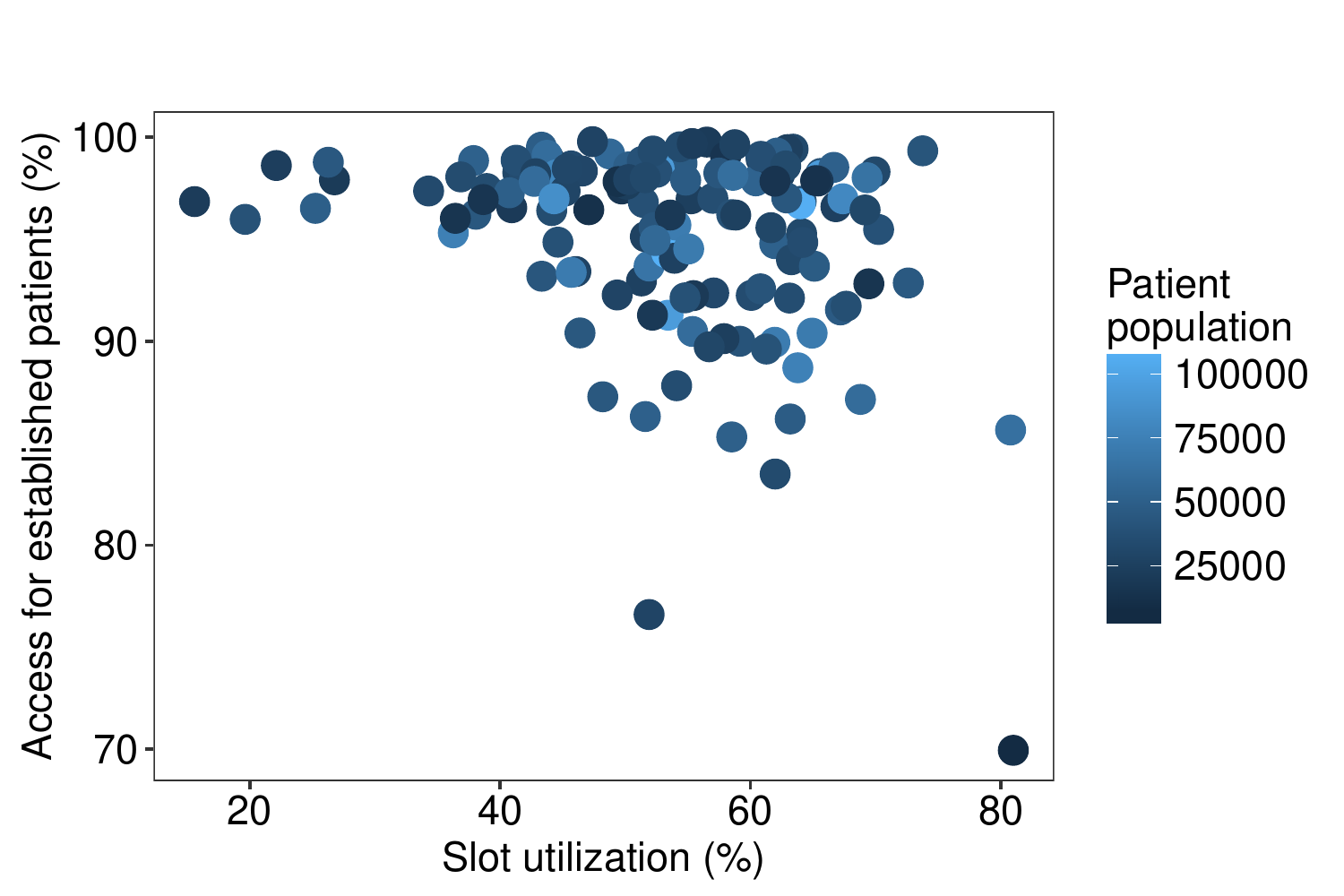}
	\caption{Access to primary care for returning patients vs. slot utilization at VA facilities in 2013. Color coding indicates the number of primary care patients cared for at each facility.}
	\label{AUEst}
\end{figure}

The American Academy of Family Physicians defines primary care as care by providers who are trained for comprehensive first contacts and continuing care for patients with any undiagnosed sign, symptom, or health concern \cite{AAFP}. Access to primary care, care quality, and health service efficiency are important dimensions of healthcare system performance \cite{Commonwealth2011primary}. One way to improve the quality of health service delivery is to establish efficient patient flow to and within healthcare facilities \cite{dixon2015patient}. 
The lead-time to an appointment, measuring the difference between a patient's desired and actual appointment dates, is known as \textit{indirect} waiting time and the waiting time at the clinic on the day of the appointment is known as \textit{direct} waiting time \cite{gupta2008appointment}. In the U.S., the average indirect waiting times for 2014 varied from five days in Dallas to 66 days in Boston \cite{Merritt2014Hawkins}. When it comes to the Department of Veterans Affairs (VA), the largest healthcare system in the U.S., access to care has been a struggle \cite{VAAccess}. Figure \ref{AUEst} depicts the relationship between access to primary care and appointment slot utilization (defined as the percentage of total available provider appointment slots that are actually used for providing care) at VA facilities across the nation. The VA defines access with a binary measure that indicates whether a returning patient has been given an appointment within 14 days of the desired appointment date. It is clear from the plot that access across facilities varies widely in spite of appointment slot utilization being less than 60\% for the vast majority of the facilities. 

Access to primary care is expected to improve patient health outcomes, reduce overall healthcare costs, and increase health equality between population groups \cite{Access2014NACHC}. An analysis by Prentice et al. \cite{prentice2007delayed} of facility-level data from 89 VA medical centers merged with patient-level data from geriatric outpatient clinics revealed that long access delays have a significant impact on negative health outcomes such as mortality. 
In addition, appointment delays can lead to attrition in the number of patients using a facility, and they present a lost opportunity to treat patients on time \cite{pomerantz2008improving}.

Wellstood et al. \cite{wellstood2006reasonable} report that direct waiting time in primary care clinics is a significant barrier to access and continuity of care. A survey by Software Advice \cite{Patient2013Advice} shows that more than 40\% of patients are willing to visit another provider, compromising continuity of care, in order to have shorter direct waiting times. The study also shows that while 45\% of patients are able to see their provider within 15 minutes upon arrival, 15\% of patients wait more than 30 minutes. Another study by Anderson et al. \cite{anderson2007willing} shows that around 25\% of patients wait more than 30 minutes in a clinic to visit their primary care provider, and the study indicates that overall satisfaction is inversely proportional to waiting time.

The Institute of Medicine considers mismatches between supply and demand to be one of the main causes of delays in access to healthcare \cite{kaplan2015transforming}. While demand for healthcare in the U.S. is expected to increase by 29\% between 2005 and 2025 due to population growth and aging, it is estimated that the number of adult primary care practitioners will only increase by 2$\sim$7\% during the same period \cite{bodenheimer2010primary}. Balancing supply and demand in the healthcare environment is usually done through appointment scheduling systems, and these studies indicate that factors such as a provider's typical service time per patient type and patient preferences regarding the day and time of their appointment, patient tardiness, and patient no-shows are uncertainties in the appointment scheduling problem. Ignoring these may result in scheduling rules that are sub-optimal or infeasible in real clinical settings.

The goal of this study is to bridge the gap between appointment scheduling and patient flow in primary care clinics. We propose a two-stage stochastic programming model to develop robust scheduling templates that account for uncertainties in demand volume, patient preferences for appointment times, provider schedule, overtime (i.e., outside-of-regular-hours work) restrictions, and patient no-show rates to balance indirect and direct waiting times. The two-stage model uses the sample average approximation (SAA) method for computational efficiency, which asymptotically converges to an optimal solution. The resulting template allocates the expected demand for different days and appointment slots based on patient types and resource availability during the booking horizon.

\begin{figure}[!t]
	\centering
	\includegraphics[scale = 0.35]{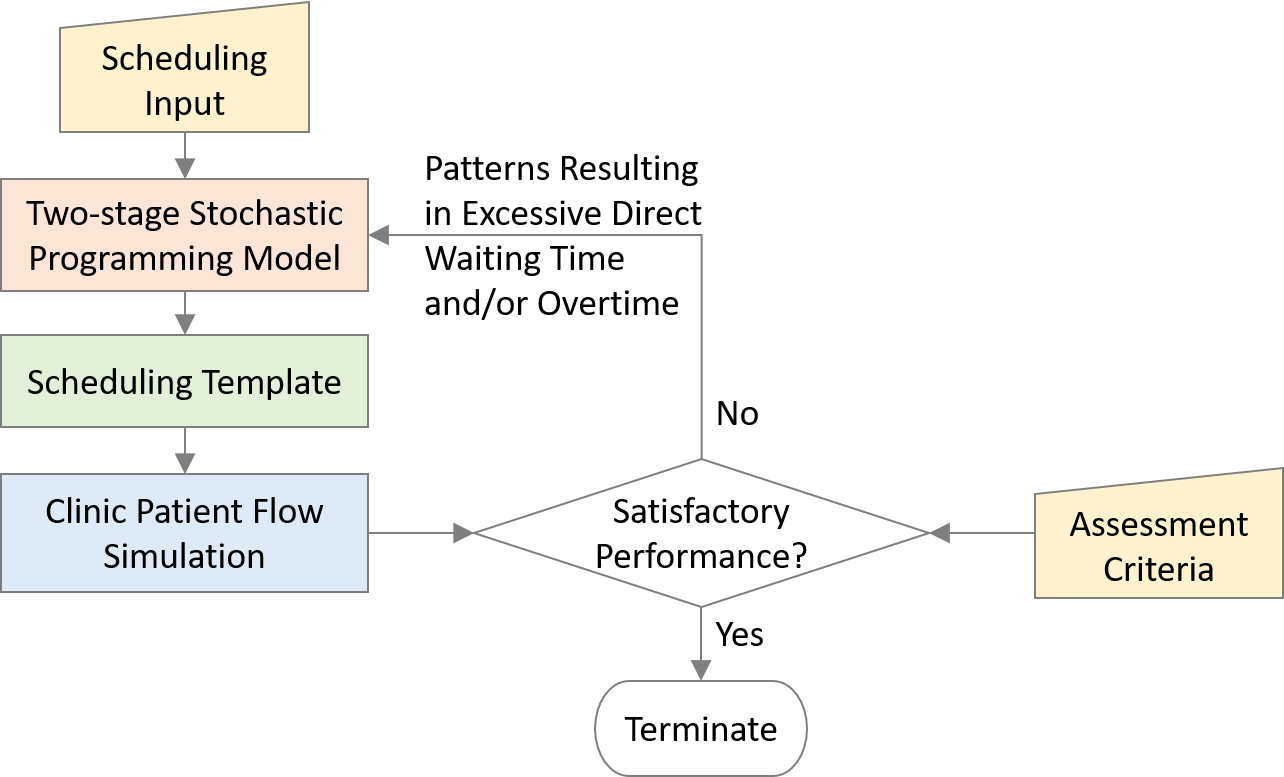}
	\caption{Proposed approach for producing robust appointment scheduling templates}
	\label{Overall}
\end{figure}

Figure \ref{Overall} shows the components of the proposed approach. First, using provided input, the two-stage stochastic programming model produces a scheduling template that minimizes indirect waiting time. The clinic patient flow simulation model then evaluates the template for satisfactory direct waiting timer performance and reports any template sequences that cause direct waiting times and overtime that exceeds specified thresholds (assessment criteria). The model is iterated with new constraints to avoid the violating sequences, and eventually, a template is obtained that minimizes indirect waiting time while balancing direct waiting times. The optimal scheduling template can then be used by the scheduling call centers to manage appointments.

The major contributions of this study are as follows: 1) We propose an integrated two-stage stochastic programming and simulation approach to developing robust scheduling templates for primary care clinics that minimize indirect and direct waiting times. 2) We use simulation modeling to manage patient flow in the clinic by introducing sequencing rules that control patients' waiting time and the provider's overtime. 3) We provide an index policy for appointment scheduling in the call center that considers several factors, such as patients' preferences for the day and time of their appointment and patient and clinic appointment cancellations. 4) We validate the proposed methods using data from real-world clinics and corroborate the efficiency of our proposed model compared to existing approaches in the literature. 

The rest of this paper is structured as follows. Section \ref{Literature review} provides a review of related literature. Section \ref{Problem description} describes the assumptions and the uncertainties that are considered in our model. This section also presents the model formulation and our simulation approaches to clinic patient flow and call center scheduling. Section \ref{Numerical study} defines various performance measures for appointment scheduling and derives practical guidelines through a numerical study of a VA primary care clinic. Finally, we summarize our conclusions and discuss directions for future work in Section \ref{Conclusion}.

\section{Literature review} \label{Literature review}

There is a rich body of healthcare operations management literature on outpatient appointment scheduling. However, prior research has mostly focused on proposing appointment scheduling systems to manage patient flow within clinics (i.e., direct waiting time) but not to effectively balance both direct and indirect waiting times.   

\subsection{Clinic patient flow}
\subsubsection{Patient flow measures}

Muthuraman et al. \cite{muthuraman2008stochastic} proposed a stochastic overbooking model to optimize appointment scheduling in an outpatient clinic where patients have different no-show probabilities. Their objective function captures patient waiting time, provider overtime (i.e., work outside-of-regular-hours), and idle time. Zeng et al. \cite{zeng2010clinic} maximized clinics' expected profit based on revenue from patients and the costs of patient waiting time, provider overtime, and idle time, with patients having different no-show probabilities. The authors observed that the performance of scheduling practices using homogeneous overbooking models based on the mean value of show-up probabilities is not good enough. Chakraborty et al. \cite{chakraborty2013sequential} developed a sequential scheduling algorithm to minimize the total expected cost resulting from patients' waiting time and providers' overtime using stochastic service times. They showed that their model leads to higher profits and less overtime than policies that consider service periods to be pre-divided into slots.

\subsubsection{No-shows and overbooking}

Patient no-shows are a major challenge in outpatient clinics. To mitigate the negative impact of no-shows on scheduling practice, Laganga et al. \cite{laganga2012appointment} developed an appointment scheduling approach using overbooking to balance patients' waiting time and providers' overtime. They concluded that it is impossible to draw general conclusions about constructing overbooking schedules. Zacharias et al. \cite{zacharias2014appointment} proposed an overbooking model to mitigate the negative impact of patient no-shows on clinic performance when patients have different no-show probabilities. The authors studied static and dynamic scheduling problems and showed that patients' heterogeneity in no-show rates has a large negative impact on the scheduling process. 

\subsection{Indirect waiting time}
\subsubsection{Advanced access scheduling}

Clinics tend to use \textit{advanced access} systems to reduce patients' indirect waiting time. In these systems, patients are given appointments on or near their desired date. Taking into account patients' no-show and appointment cancellation behavior, Liu et al. \cite{liu2010dynamic} proposed a dynamic scheduling policy for an outpatient clinic and showed that an advanced access scheduling policy performs better when the demand rate is relatively low. Dobson et al. \cite{dobson2011reserving} examined the effect of keeping some slots open for same-day demand in primary care clinics on two quality measures: the average number of same-day demands that are not served during normal working hours and the average number of non-urgent patients in the queue. They demonstrated that encouraging non-urgent patients to call for same-day appointments is an important factor when implementing advanced access systems in primary care clinics. Qu et al. \cite{qu2012mean} derived the selection percentage for open appointments in an advanced access system by using a mean-variance approach. Their results indicated that when both the demand rate and the no-show rate are high for appointments that are reserved for routine patients, there are one or more Pareto optimal percentages of open appointments that can decrease the variability in the number of patients seen. 

\subsubsection{Patient choice}

Patient scheduling choices can impact appointment delays. Gupta et al. \cite{gupta2008revenue} developed a Markov decision process to manage access to care when patients can choose between accepting a same-day or a future appointment. The authors provided optimal solutions for clinics with single and multiple providers. Wang et al. \cite{wang2011adaptive} studied clinic revenue optimization by finding the optimal balance between the number of slots that should be kept open for same-day demand and the number of slots for routine patients while considering preferences regarding providers and appointment times. Their model is limited to one day and hence does not consider the interactions between multiple scheduling days.

\subsection{Summary}
Our work is closer to Luo et al.'s \cite{luo2015tandem} research in which they developed a tandem queue model to study the relationship between the appointment queue (indirect waiting time) and the service queue (direct waiting time). The main research question that we address is: How can primary care practices schedule patients to ensure that the patients experience minimal delays in getting their appointments while making the patient flow in the clinic as smooth as possible? Our work is different from the studies discussed above in several important respects. First, we consider the indirect waiting time of patients who may call in advance to book an appointment over a planning horizon $T$. Second, in the optimization model, we consider three patient flow measures: patients' direct waiting time and the provider's overtime work during lunchtime and after regular hours. Third, we account for patients' preference for appointment dates and times. Finally, we focus on scheduling \textit{templates} that are easy to employ by appointment call center staff given their \textit{stability} for extended periods (e.g., a month or a quarter) and also allowing providers the benefit of a stable work pattern with considering their clinical and administrative practice. Like other studies, we also consider the effect of a patient's no-show behavior on patient scheduling and patient flow measures. 

\section{Problem description and model formulation} \label{Problem description}

In this section, we provide an overview of the problem setup and our assumptions. We also present the notation used in the model, followed by the proposed risk-neutral two-stage stochastic programming model for developing appointment scheduling templates.

\subsection{Problem setup and assumptions} \label{setup}
\begin{figure*}[!t]
	\centering
	\includegraphics[scale = 0.5,trim={0 0 0 1cm},clip]{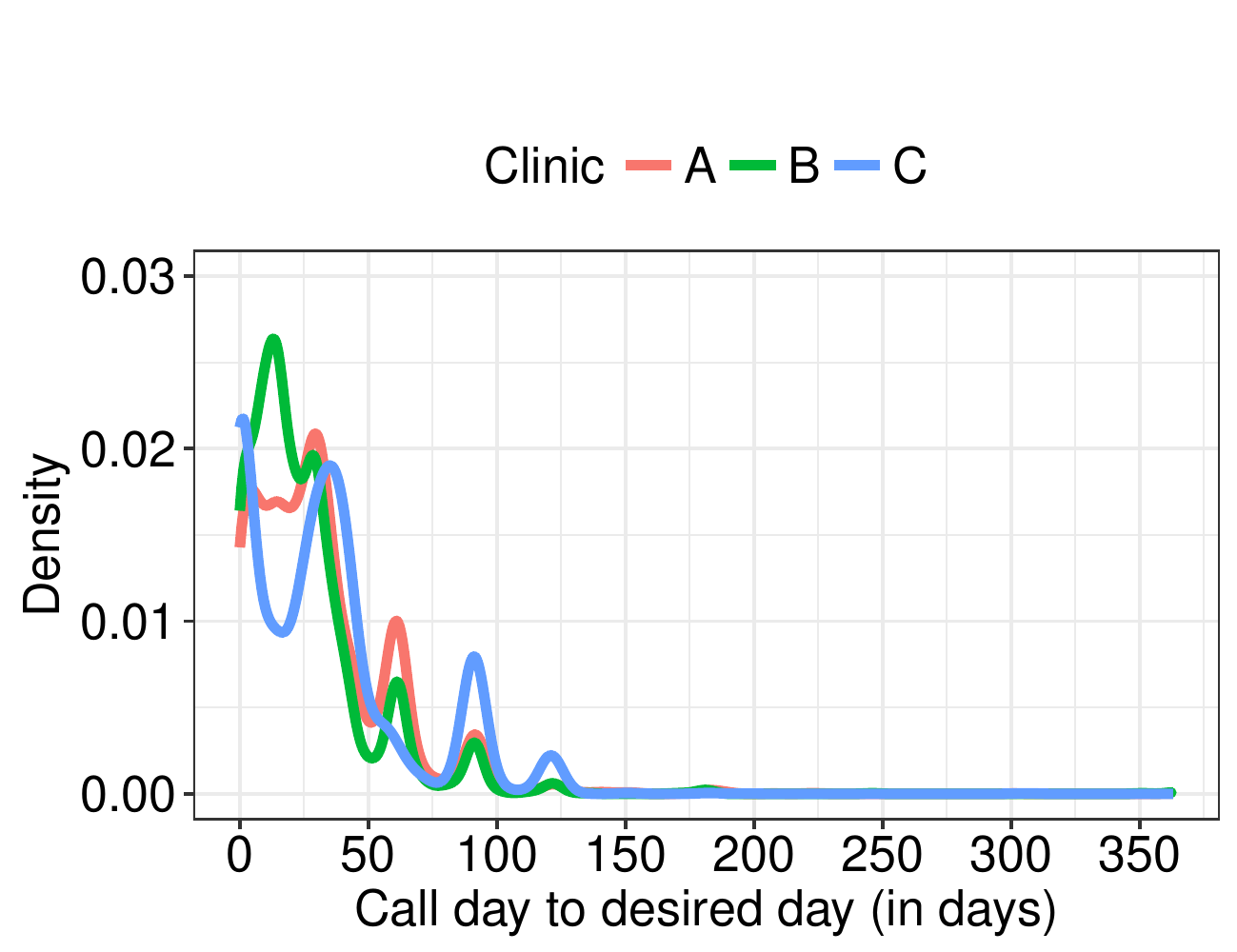}
	\caption{Time (in days) between call date and desired appointment date for patients in three VA primary care clinics. The clinic names are coded.}
	\label{Desired}
\end{figure*}

We study primary care outpatient clinics with single providers.  Although solo practice is becoming less popular, more than half of family physicians still work in solo and small practices \cite{liaw2016solo}. Without loss of generality, we assume that there are five working days per week (Monday through Friday) and each working day has eight provider working hours, made up eight 60-minute appointment slots. Each day is divided into morning and evening sessions, each of which lasts four hours (8AM to Noon and 1 to 5PM), with a one hour lunch break. Without loss of generality, we assume that there are 20 working days in a month, resulting in 40 sessions/month. Patients call to book appointments with a preference for the appointment date (future or same-day) and time (early morning or closer to lunch hour or late afternoon). While much of the academic literature assumes that patients call for an appointment on the day when they want the appointment, data shows that most primary care patients often call in advance to make future appointments. Figure \ref{Desired} shows how early patients call to request appointments in three different VA primary care clinics in the Midwest. 

As is typical with primary care, we allow multiple patient types (e.g., new patients vs established patients, annual physicals, etc.) and we allow different service times for nurses and providers based on the patient type. Clinics are allowed to cancel appointments (e.g., due to lab result delays or a provider's absence). These cancellations need to be managed since they will increase patients' dissatisfaction and the staff's future workload. Patients may also cancel their appointments or not show up at all for their visit. We allow overbooking as a means to compensate for patients' no-show behavior with respect to providers' appointment slot utilization. However, excessive overbooking increases patients' direct waiting time. 

Without loss of generality, we assume that patients' preferences cannot be denied (e.g., appointment time). Otherwise, patients will seek care in a specialty care clinic or an emergency department, both of which are more costly than primary care. Even though providers can work overtime, the mathematical model we develop needs to moderate the effect of capacity shortage due to other responsibilities the provider may have. Also, patients typically have different preferences for their appointment dates, and there are uncertainties regarding the number of patients who will call each day during the planning horizon and their desired appointment dates.

The two-stage model we develop minimizes patients' indirect waiting time while considering patient flow within the clinic. Also, the model ensures that the provider is not overloaded with excessive cumulative workloads in any morning or evening session by tracking expected overtime during lunch hour and work past the end of the work day. Finally, we consider a finite rolling scheduling horizon and seek weekly scheduling templates for they are typical in practice. We allow the weekly templates to vary from month to month to account for any seasonal differences in demand patterns.

\subsection{Two-stage stochastic programming model} \label{model}

Stochastic programming is a branch of optimization that assumes that some of the model parameters and coefficients are unknown and that only their probability distribution can be estimated. The most widely used stochastic programming model is two-stage stochastic programming. In this model, the first-stage decision variables are ``here-and-now'' decisions that are determined before observing the realization of uncertainties, and the second-stage decision variables are selected after exposing the first-stage variables to the uncertainties. The goal is to determine the values for first-stage decisions in a way that minimizes the first- and second-stage objective function values. 

We consider uncertainties in the number of patients seeking care, the days when they call, their desired appointment dates and times, and their no-show rates in the model. The first-stage decisions in our proposed two-stage model determine the number of patients of each patient type allowed in each appointment slot and session based on the provider's maximum tolerable cumulative patient \textit{complexity} (determined in terms of expected cumulative service times). Based on the first-stage patient allocation decisions and the realization of uncertainties in the second stage, appointment scheduling decisions are made in the second stage that minimize patients' total indirect waiting time. The output of the two-stage stochastic programming is a weekly scheduling template for the booking horizon. 

\subsection{Model notation}

We consider a set of $R$ patient types, indexed by $r$, each with an average complexity $c_r$ and no-show probability $p_r$, who phone to request an appointment. The planning horizon has multiple working days, denoted by $D$ and indexed by $d$. Each day has two sessions, denoted by $S$ and indexed by $s$. Within each session, there are multiple appointment slots, denoted by $A$ and indexed by $a$. In order to manage the patient flow in the clinic and handle different clinical tasks, maximum patient complexities are considered for each appointment slot and session, respectively denoted as $\kappa$ and $\eta$. The numbers of patients of type $r$ that can be scheduled in each appointment are given by the set $L(r)$, indexed by $l$. 

For each patient type $r \in R$, let $l \in L(r)$; then the parameter $m_{r,l}$ denotes the discrete number of possible patients of type $r$ that can be scheduled in each appointment slot. To help maintain a rolling planning horizon for the weekly scheduling grid template, the parameter $\xi_{r,a}$ denotes the number of previously booked patients in the template. Let $\tilde{\omega}$ be a random variable representing the uncertainties in the two-stage model, and let $\omega$ be a realization of $\tilde{\omega}$. The first-stage decision variables $x_{r,a}$ and $z_{r,a,l}$ respectively determine the number of patients of type $r$ that can be scheduled in appointment slot $a$ and whether $l$ patients of type $r$ can be scheduled in appointment slot $a$. The second-stage decision variable $y_{r,d,d'}$ assigns patients of type $r$ that call on day $d$ and request an appointment on day $d'$. Table \ref{notations} summarizes the notation that is used for the two-stage stochastic programming model. 

\begin{table}[H]
	\caption{Model Notation} 
	\begin{tabular}{p{0.2\linewidth}p{0.75\linewidth}}
		\hline
		\textbf{Symbol} & \textbf{Description}\\
		\hline
		Sets:\\
		$R$			& Set of patient types, indexed by $r \in R$ \\
		$A$			& Set of appointment slots, indexed by $a \in A$ \\
		$S$		    & Set of sessions, indexed by $s \in S$ \\
		$G$		    & Set of template sequences in which patient flow constraints are not met for direct waiting time and provider overtime work thresholds, indexed by $g \in G$ \\
		$D$	        & Set of days, indexed by $d \in D$ \\
		$L_{r}$     & Set of numbers of patients of type $r$ that can be scheduled in each appointment slot, indexed by $l \in L_{r}$ \\
		$\Omega$    & Set of scenarios, indexed by $\omega \in \Omega$ \\
		\hline
		Model Parameters:\\
		$c_r$       & Average complexity of patient type $r$ \\
		$\kappa$    & Maximum acceptable cumulative patient complexity for each appointment slot \\
		$\eta$      & Maximum acceptable cumulative patient complexity for each session \\
		$p_r$       & Average no-show probability for patient type $r$ \\
		$m_{r,l}$   & Number of patients of type $r$, $l \in L_{r}$ \\
		$\xi_{r,a}$ & Number of scheduled patients of type $r$ in appointment slot $a$ \\
		$f_{r,d}(\omega)$ & Number of patients of type $r$ who asked for an appointment on day $d$ in scenario $\omega$ \\
		$\epsilon$ & User parameter \\
		\hline
		First-stage Variables: \\
		$x_{r,a}$   & Number of patients of type $r$ who can be scheduled in appointment slot $a$ \\
		$z_{r,a,l}$ & 1 if $l$ patients of type $r$ can be scheduled in appointment slot $a$; 0 otherwise \\
		\hline
		Second-stage  Variables: \\
		$y_{r,d,d'}(\omega)$ & Proportion of patients of type $r$ who asked for an appointment on day $d$ and are scheduled for day $d'$ in scenario $\omega$ \\
		\hline
		\label{notations}
	\end{tabular}
\end{table}
\vspace{-3em}

The first-stage problem is represented as follows:

\begin{align}
\text{ Min } & f(x) = \mathbb{E}[\varphi(x, \widetilde{\omega})] & \label{eq1} \\
\text{s.t.} & \sum_{r \in R} c_rx_{r,a} \leq \kappa \hspace{2.6cm} \forall  a \in A, \label{eq2} \\
& \sum_{r \in R} \sum_{a \in s}  c_rx_{r,a} \leq \eta \hspace{2.1cm} \forall s \in S, \label{eq3} \\
& x_{r,a} \geq \xi_{r,a} \text{ } \hspace{2.1cm} \forall r \in R, a \in A, \label{eq4} \\ 
& x_{r,a} = \sum_{l \in L_r} m_{r,l}z_{r,a,l} \hspace{0.7cm} \forall r \in R, a \in A, \label{eq5} \\
& \sum_{l \in L_{r}} z_{r,a,l} = 1 \hspace{1.7cm} \forall r \in R, a \in A, \label{eq6}\\
& \sum_{r \in R,a \in A,l \in L_{r} \subseteq G} z_{r,a,l} \leq |G|-1, &  \label{eq7}
\end{align}

\begin{align}
&  x_{r,a} \in \mathbb{Z}^+ \hspace{3cm} \forall r \in R, a \in A, & \nonumber\\
&  z_{r,a,l} \in \{0,1\} \hspace{1.4cm} \forall  r \in R, a \in A, l \in L_{r}. \label{eq8}
\end{align}

For a given first-stage solution $x$ and the realization $\omega$ of random variables, the second-stage recourse function is as follows:

\begin{align} 
\text{ Min } & \varphi(x,w) =  & \nonumber\\
& \sum_{\substack{r \in R \\ d,d' \in Dd \leq d'}}   w_ry_{r,d,d'}(\omega)f_{r,d}(\omega)[(d'-d)^{(1+\epsilon)}] & & \label{eq9}\\
\text{s.t.} & \sum_{d \in D: d \leq d'} (1-p_r)y_{r,d,d'}(\omega)f_{r,d}(\omega) \leq \sum_{a \in d'} x_{r,a}(\omega) & \nonumber\\
& \hspace{3.4cm} \forall r \in R, d' \in D, \label{eq10} \\
& \sum_{d' \in D: d \leq d'} y_{r,d,d'}(\omega) = 1 & \nonumber\\ 
& \hspace{3.4cm} \forall r \in R, d \in D, \label{eq11}\\
&  0 \leq y_{r,d,d'}(\omega) \leq 1 & \nonumber \\ 
& \hspace{1.7cm} \forall r \in R, d,d' \in D: d \leq d' \label{eq12}.  
\end{align}

The objective function minimizes the expected indirect waiting time for patients. The difference between the desired and actual appointment dates for a patient is penalized using a super-linear function in order to favor ``fairness'' in assigning lengths of delay to the patients.    

A primary care provider's threshold in terms of the cumulative patient complexity that can be handled in each appointment slot and each scheduling session is represented in constraints \eqref{eq2} and \eqref{eq3}, respectively. Since this is a rolling planning horizon problem, constraints \eqref{eq4} fill the slots based on a commitment to previously scheduled appointments. Constraints \eqref{eq5}, \eqref{eq6}, and \eqref{eq7} are sequencing rules to address patient flow in the clinic, and they are added dynamically based on recommendations from the clinic patient flow simulation model. Constraints \eqref{eq5} and \eqref{eq6} determine the maximum number of patients of each type that can be scheduled in each appointment slot, and constraints \eqref{eq7} ensure that sequences that have violated patient flow thresholds do not occur in the scheduling template. Constraints \eqref{eq8} are integer and binary value constraints for the first-stage variables. Constraints \eqref{eq10} ensure that patient appointments are provided based on the scheduling template resulting from the first-stage model. Constraints \eqref{eq11} ensure that no patient request is denied, and constraints \eqref{eq12} confirm that the second-stage variables are proportion values between 0 and 1. 

\begin{figure*}[!t]
	\centering
	\includegraphics[scale = 0.6]{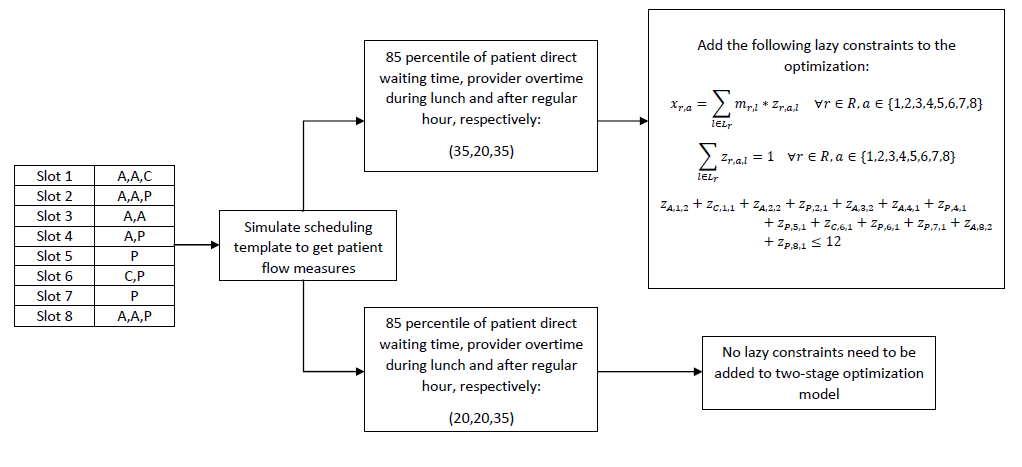}
	\caption{An example of adding lazy constraints to the two-stage optimization model based on the clinic patient flow simulation performance}
	\label{Seqplot}
\end{figure*}

Figure \ref{Seqplot} represents an example of adding ``lazy constraints'' to the two-stage stochastic model. Based on the scheduling template from the two-stage model, the patient flow simulation determines the percentiles (85th percentile in this example) of patients' direct waiting time and the provider's overtime work during lunchtime and after regular hours, and if any of these violate the predetermined thresholds, which are 30 minutes for patients' expected direct waiting time, 45 minutes for the provider's expected overtime work during lunch, and 60 minutes for the provider's expected overtime work after regular hours, a lazy constraint is added to the first stage of the optimization model to eliminate such sequences. As shown in Figure \ref{Seqplot}, when the simulated performance values are 20, 20, and 35 for patients' direct waiting time and the provider's overtime during lunch and after regular hours, respectively, none of the thresholds is violated and so lazy constraints are not added to the two-stage stochastic model. However, when these values are 35, 20, and 35, respectively, the threshold for patients' direct waiting time is violated, and the corresponding lazy constraint is added to the optimization model. Here, $x_{r,a}$ is the number of patients of type $r$ who can be scheduled in appointment slot $a$, and $z_{r,a,l}$ is equal to 1 if $l$ patients of type $r$ can be scheduled in appointment slot $a$. Acute, chronic, and preventive patient types are represented as `A', `C', and `P', respectively. Whenever a lazy constraint is added to the first stage, the model is re-optimized.

\vspace{-6mm}

\subsection{Clinic patient flow simulation}

The patient flow simulation for the clinic is executed for each day based on the scheduling template proposed by the optimization model described in the previous section. Patients' direct waiting time and the provider's overtime work during lunch time and after regular hours are measured. The clinic patient flow involves two stages: time with the nurse and time with the provider. Once a patient walks into the clinic, he or she waits in the lobby for the nurse to become available. After being visited by the nurse, the patient waits for the provider in the exam room. We assume that patients who are scheduled for a particular day must be served before the end of that day even if the provider has to work overtime. Patients are assumed not to leave the exam room until the provider finishes all required tasks. Also, patients may arrive late for their appointments. By arriving late, patients may increase the waiting times for the patients that follow them. Therefore, we assume that patients are called in the order of their arrival time. 

Algorithm \ref{patientflow1} presents pseudocode for our approach based on the two-stage optimization model and simulation. In step 1, the two-stage stochastic programming model is solved, and a scheduling template is obtained from the first-stage solution. Using the scheduling template, each day in the planning horizon is simulated under given input with $m=200$ replications in step 2, and specific percentiles of patient flow measures are calculated in step 3. If any of the measures violate the patients' or the provider's thresholds, a new set of constraints \eqref{eq5}, \eqref{eq6}, and \eqref{eq7} are added to the model, which is then re-optimized. The process repeats until the system reaches a state where the indirect waiting time is minimized without violating any patient flow thresholds.  

\begin{algorithm}[H]
	\SetAlgoLined
	\textbf{Step 1:} Optimize the two-stage stochastic model \eqref{eq1} - \eqref{eq4}, \eqref{eq8} - \eqref{eq12}; \\
	\textbf{Step 2:} Run simulation model for 200 replications for each day in the planning horizon; \\
	\textbf{Step 3:}
	\For{$i$ in $D$}{
		Calculate $\alpha^{th}$ percentile for patients' direct waiting time and provider's overtime during lunch and after regular hours;\\
		\uIf{estimates $\geq$ any of the  thresholds}{
			Construct constraints \eqref{eq5}, \eqref{eq6}, and \eqref{eq7} and add them to the two-stage stochastic model;\\ 
	}}
	\caption{Clinic patient flow simulation}
	\label{patientflow1}
\end{algorithm}

\subsection{Call center simulation}

To evaluate the efficiency of the two-stage scheduling template, we simulate the practice at the appointment call center. The call center simulation uses either the total available capacity or the scheduling template's allocation. An ``index scheduling policy'' is used to allot or cancel patients' appointments. The scheduler at the call center estimates the priority of available appointment slots that should be offered to a patient based on his or her desired date and patient type. When a patient requests an appointment for a desired date, the policy calculates the index based on the slot's remaining capacity in increasing order for each appointment slot based on the slot's proximity to the desired date. To generate a patient's choice regarding accepting an appointment, a random number from the uniform distribution $\mathcal{U}(0,1)$ is generated and compared to an ``acceptance'' threshold. If the random number is greater than this threshold, the patient accepts the corresponding appointment slot; otherwise, another slot is offered and the process continues until the patient accepts. Appointment cancellations are handled similarly: A random number is generated from $\mathcal{U}(0,1)$, and if the random number is less than the clinic's cancellation rate, the patient or clinic cancels the appointment and the patient is removed from the scheduling grid.

\section{Case study and insights} \label{Numerical study}

We used real data from a U.S. Midwest VA primary care clinic to estimate the number of weekly requests, patients' no-show probabilities, appointment cancellation probabilities, the daily distribution of patients' calls, the distribution of patients' desired appointment days, and the distribution of the time between call dates and desired dates. The data suggests that patients often call with the same probability on different weekdays, but fewer patients ask for appointments on Mondays while more ask for Fridays. Patients may ask for appointments up to four weeks in advance, but around 65\% of the patients want an appointment within one week. No-show and cancellation rates vary across months, and the average no-show and cancellation rates for this clinic are 10\% and 17\%, respectively. Other parameters are listed in Table \ref{Base}. 

\vspace{-3mm}
\begin{table}[H]
	\caption{Base problem parameters} \label{Base}
	\begin{tabular}{p{0.9\linewidth}p{0.1\linewidth}}  
		\hline
		Maximum patient complexity that the provider can handle in an appointment slot = 0.96 \\ 
		Maximum patient complexity that the provider can handle in a session = 2.8 \\
		Threshold for patients' acceptance of offered appointment slot = 0.2 \\ 
		Average complexity of different patient types = [Acute: 0.29, Chronic: 0.32, Preventative 0.36] \\ 
		Threshold for patients' direct waiting time = 30 minutes\\ 
		Threshold for spillover amount to provider's lunch time = 45 minutes \\
		Threshold for provider's overtime = 60 mintues\\ 
		Percentile of patient flow metric distributions in clinic patient flow simulation = 85\% \\
		Patient arrival time distribution = $\mathcal{N}(-16.62,27)$ \\
		Booking horizon = 60 days \\
		Planning horizon = 300 days \\
		\hline
	\end{tabular}
\end{table}
\vspace{-3mm}

We considered three different patient types (acute, chronic and preventive) based on a study by Yarnall et al. \cite{krause2009family}, who  used the National Ambulatory Medical Care Survey for 2003 to determine the visits for these patient types. Similarly, we used an empirical study by Oh et al. \cite{oh2013guidelines} to determine the amount of time the nurse spends with patients of each type. The service time with the nurse and the provider for each patient type follows a log-normal distribution, as suggested by Cayirli et al. \cite{cayirli2006designing}. Table \ref{serv} reports the service time with the nurse and the provider for each patient type along with the percentage of each patient type in the provider's panel. 

\begin{table}[!htbp]
	\begin{center}
		\caption{Expected service time with nurse and primary care provider}
		\label{serv}
		\centering
		\begin{tabular}{| c | c | c | c |} 
			\hline
			Visit & (\%) of & Nurse Time & Provider Time \\
			Type & Visits & (mins) & (mins) \\
			\hline
			Acute & 49.3 & 11.3 (8.3) & 17.3 (8.7) \\ 
			\hline
			Chronic & 36.1 & 12.6 (8.8) & 19.3 (9.2) \\
			\hline
			Preventive & 14.6 & 13.9 (11.3) & 21.4 (11.8) \\
			\hline
		\end{tabular}
	\end{center}
\end{table}

Delays in arrivals for appointments are prevalent in outpatient clinics. In this study, we used the normal distribution $\mathcal{N}(-16.62, 27.07)$ in minutes, as estimated by Cayirli et al. \cite{cayirli2006designing}. These authors collected data from a primary healthcare clinic in a New York metropolitan hospital and used the Kolmogorov--Smirnov test to estimate the parameters. A negative average indicates that on average, patients arrive earlier than the starting time for their appointment. 

We compare the performance of the appointment scheduling based on two-stage stochastic programming with two sequencing rules that have been proposed in the literature. These rules are shortest processing time (SPT) and low coefficient of variation (CV) in the beginning (LCVB). SPT schedules patients in increasing order of mean service times, while LCVB schedules in increasing order of the CV ($\sigma/\mu$) of the service time \cite{bhattacharjee2016simulation}. For the given planning horizon, appointments were assigned based on each of these two approaches along with our proposed approach, and the corresponding direct and indirect waiting times were evaluated. Table \ref{temp} shows sample ``daily" templates for the SPT and LCVB scheduling policies that repeat every day during the planning horizon.

\begin{table}[!htbp]
	\begin{center}
		\caption{Sample daily templates for heuristic appointment scheduling policies}
		\label{temp}
		\centering
		\begin{tabular}{| c | c | c |} 
			\hline
			App. &  & \\
			slot & SPT & LCVB \\
			\hline
			1 & A,A,A & C,C,C \\ 
			\hline
			2 & A,A & C,C,C \\
			\hline
			3 & A,A & A,A \\
			\hline
			4 & A,A & A,A \\
			\hline
			5 & C,C,C & A,A,A \\
			\hline
			6 & C,C,C & A,A \\
			\hline
			7 & P,P & P,P \\
			\hline
			8 & P & P \\
			\hline
		\end{tabular}
	\end{center}
\end{table}

The planning horizon is 240 working days (corresponding to an year) and all performance measures are tracked and reported once the system has reached steady state (around 60 days) to discard the transient effects at the beginning from model initialization. For the two-stage model, we used the SAA method to estimate the number of scenarios that were required for representing uncertainties. SAA is a Monte Carlo simulation-based sampling procedure that approximates the expected value of the objective function by using a finite sample of scenarios \cite{mak1999monte}. Due to the rolling horizon, we used SAA on various days to find the most reliable number of scenarios. Based on the SAA results, we used 10 scenarios for our computational experiments since the gap between the upper and lower bounds was within 5\%. All computational studies were implemented using Python, and Gurobi 6.5 was used as the mixed-integer programming solver on a computer running Windows 7 with 2.6 GHz processing speed and 80 GB of RAM.

The computational study was conducted in four parts, which are discussed in the following subsections. In the first part, the trade-off between indirect and direct waiting times of patients in the outpatient clinic was evaluated. In the second part, we analyze how the indirect waiting time changes if patients are more sensitive to appointment delays, as well as the subsequent impact on show-up probabilities. In the third part, the influence of the provider type on indirect and direct waiting times is estimated by considering different provider capacities. In the last part, we study the relationship between the perishability of appointment slots in a clinic and its impact on different approaches to minimizing patients' indirect waiting time.

\subsection{Trade-off between indirect waiting time and patient flow}

\begin{figure}[!htbp]
	\centering
	\includegraphics[scale = 0.3]{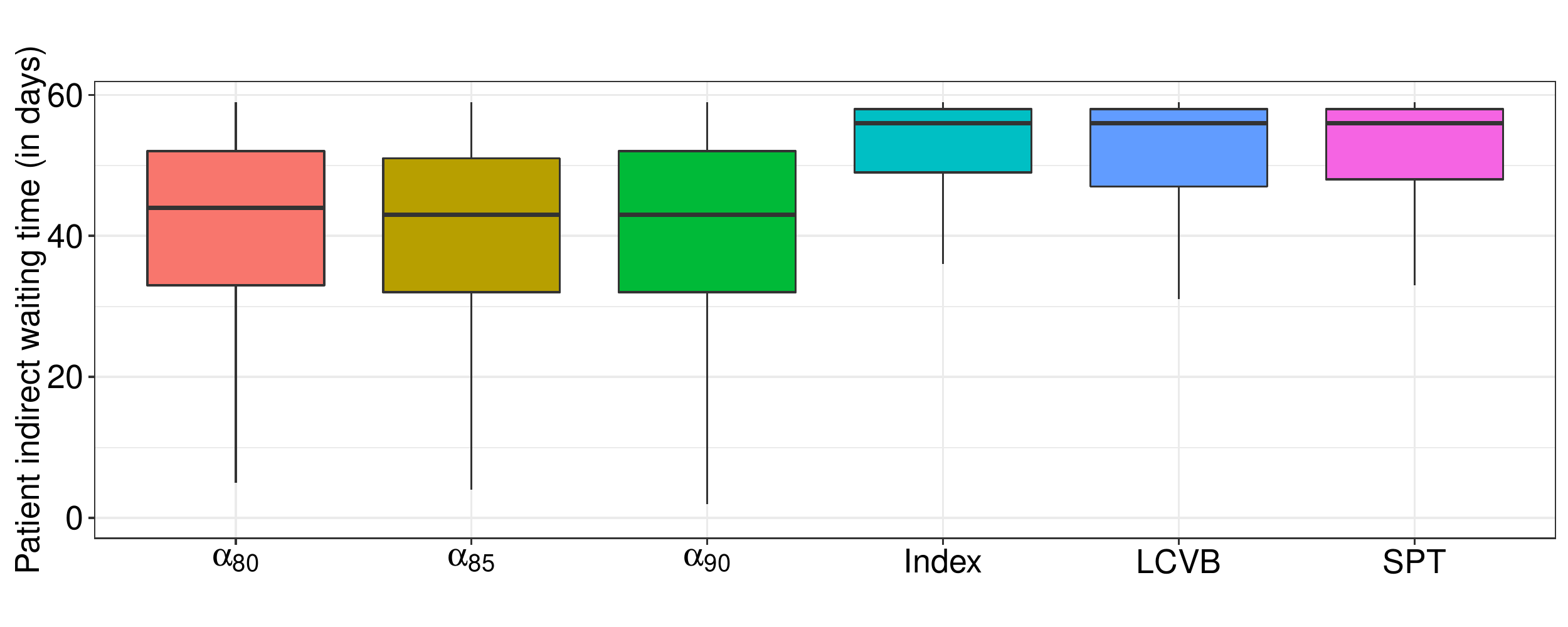}
	\caption{Patient indirect waiting time distribution under the different appointment scheduling policies}
	\label{IndirectWaitingDirect}
\end{figure}

We consider three different quantiles $\alpha$---the 80th, 85th, and 90th percentile---for the patient flow metric distributions. Figure \ref{IndirectWaitingDirect} compares the indirect waiting time distributions using our two-stage stochastic programming approach vs using baseline sequencing rules from the literature. The higher the value of $\alpha$, the more concerned the clinic manager is about patient flow in the clinic, as more patients will have waiting times that are less than their direct waiting time threshold.  

\begin{figure}[!htbp]
	\centering
	\includegraphics[scale = 0.3]{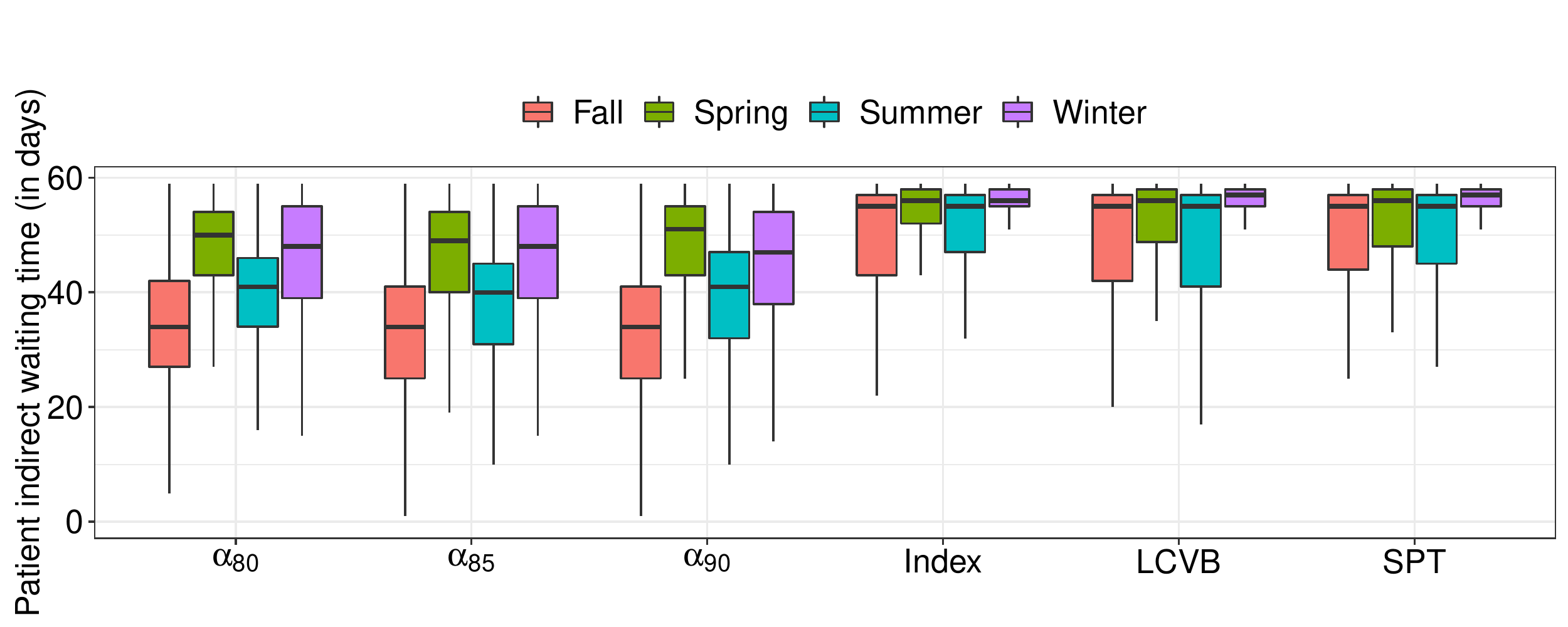}
	\caption{Patient indirect waiting time distribution under the different appointment scheduling policies in each season}
	\label{IndirectWaitingDirectSeason}
\end{figure}

\begin{figure}[!htbp]
	\centering
	\includegraphics[scale = 0.37]{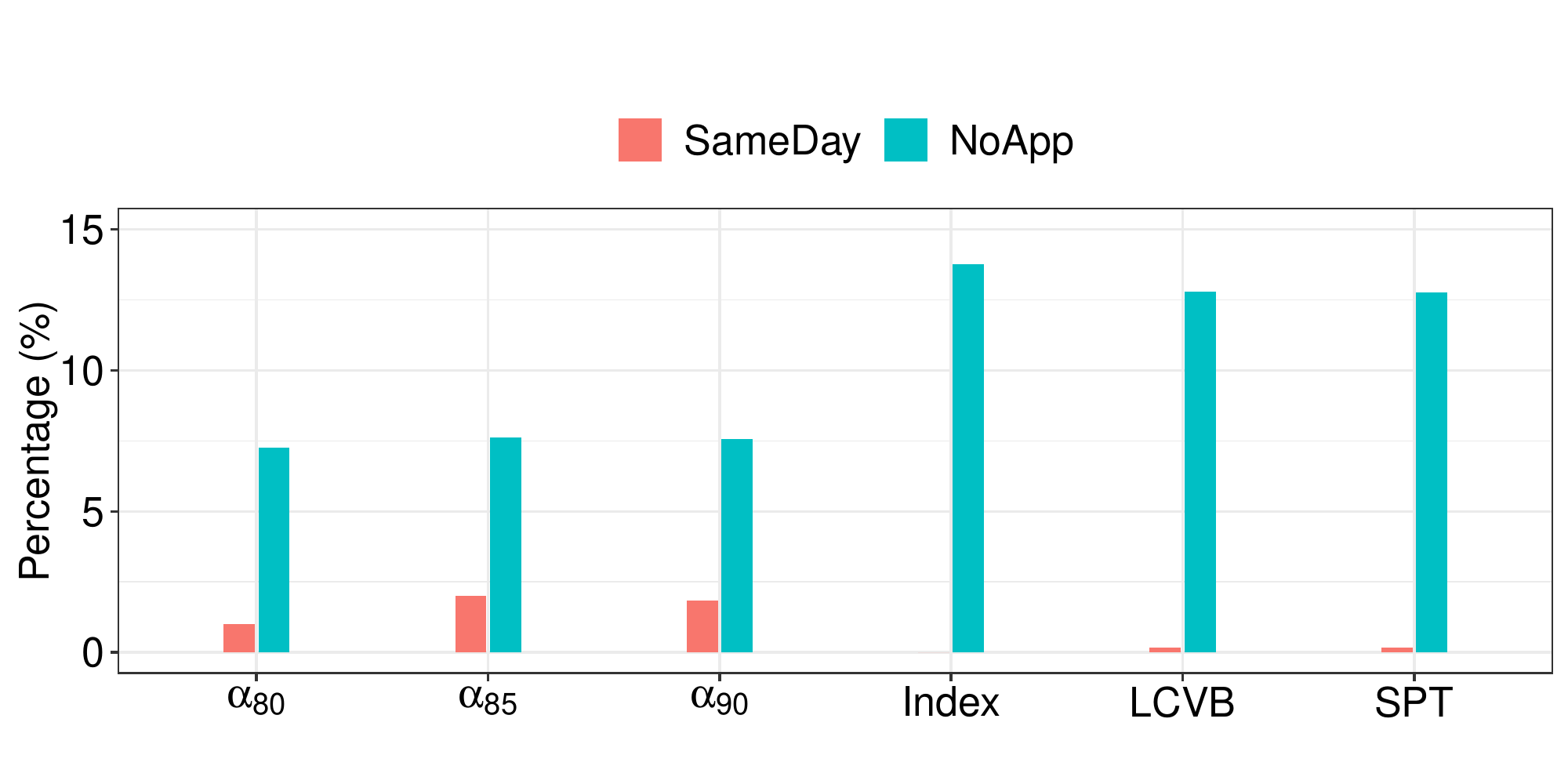}
	\caption{Percentage of patients who are given same day appointments (red) and  percentage of patients who are not given any appointment (blue) under the different scheduling policies}
	\label{IndirectWaitingDirectSameDay}
\end{figure}

The better performance of the optimal scheduling template from the two-stage stochastic programming model compared to the heuristic rules in terms of patients' indirect waiting time is shown in Figures \ref{IndirectWaitingDirect} and \ref{IndirectWaitingDirectSeason}. As Figure \ref{IndirectWaitingDirectSameDay} shows, although the difference between the optimal appointment scheduling and the heuristic policies is not significant in terms of the percentage of patients who are given same-day appointments (on average 1.6\% vs 0.1\%), the optimal policy performs better with respect to the percentage of patients who are not given any appointment with their provider (on average 7.47\% vs 13.11\%). Table \ref{temp2} represents a sample ``weekly'' scheduling template proposed by the two-stage stochastic approach for two different months across the planning horizon.   

\begin{table*}[!htbp]
	\begin{center}
		\caption{Weekly template produced by the two-stage stochastic policy in two different months across the planning horizon for the base problem.}
		\label{temp2}
		\centering
		\begin{tabular}{| c | c c c c c |} 
			\hline
			App. & & & June & & \\
			slot & Mon & Tue & Wed & Thu & Fri \\
			\hline
			1 & C,C & A,A,P & A,A,C & A,P & P,P \\ 
			\hline
			2 & A,A & No App & A,A,P & A,P & A,P \\
			\hline
			3 & A,A & A,A,C & A,A,P & A,A & A,A,P \\
			\hline
			4 & A,C,C & A,A,P & No App & A,A,A & C \\
			\hline
			5 & A,P & C,C,C & P & A,A,A & No App \\
			\hline
			6 & A,C & A,C & C,P & No App & No App \\
			\hline
			7 & A,P & A,P & C,P & A,A & No App \\
			\hline
			8 & A,A,A & A,A & A,A,P & C,C,C & No App \\
			\hline
			& & & July & & \\
			\hline
			1 & A,A,C & No App & C,P & C,P & No App \\ 
			\hline
			2 & A,C & C,P & C,C,C & A,A,A & No App \\
			\hline
			3 & A & A,A,P & No App & A,P & No App \\
			\hline
			4 & A,C,C & A,C,C & A,A,P & A,A & No App \\
			\hline
			5 & C,C,C & A,P & A,P & No App & A,A,P \\
			\hline
			6 & C,C,C & C,P & A,P & A,A,P & A,A,P \\
			\hline
			7 & C,C & C & C,C,C & A,A,P & No App \\
			\hline
			8 & No App & A,C,C & P & A,A,C & A,P \\
			\hline
		\end{tabular}
	\end{center}
\end{table*}

\begin{figure}[!htbp]
	\centering
	\includegraphics[scale = 0.55]{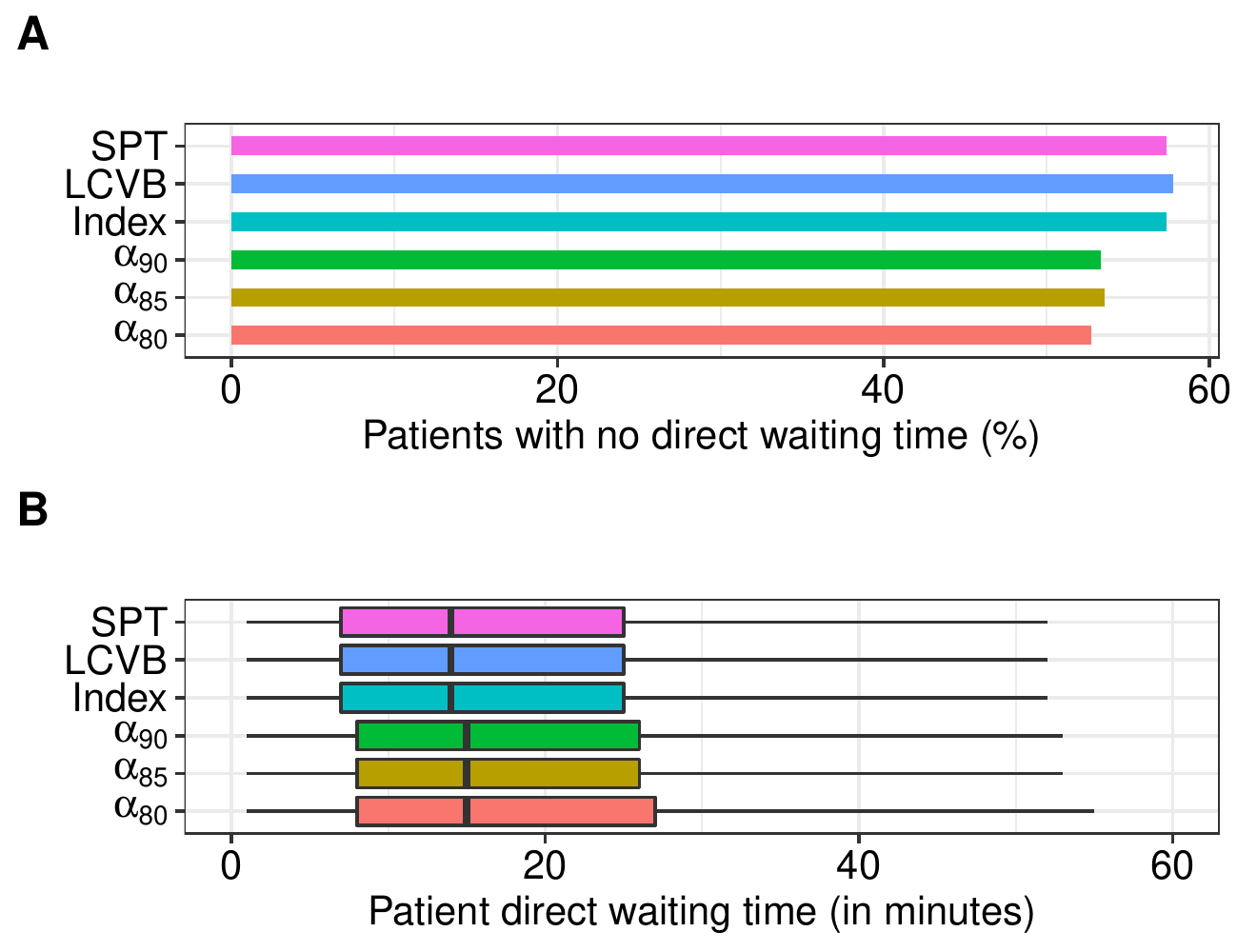}
	\caption{(A) Percentage of patients not experiencing any direct waiting time in the clinic under different scheduling policies. (B) Patients' direct waiting time distribution for patients who experience positive direct waiting times under different scheduling policies.}
	\label{DirectWaiting}
\end{figure}

On the other hand, if more patients are scheduled on a given day, then the indirect waiting time decreases but the direct waiting time increases, as shown in Figure \ref{DirectWaiting}. The major advantage of the two-stage programming approach is that it reduces the indirect waiting time without violating the threshold for the direct waiting time. Moreover, while the two-stage programming approach allows more patients to be scheduled during the last slot in each session and causes the provider to work overtime during lunch and after regular hours, the amount of extra work does not violate the corresponding thresholds. On average, the provider has to work an additional 9.25 and 4.53 minutes during lunch time and after regular hours, respectively, under the optimal scheduling policy, compared to 8.53 and 5.66 minutes under the heuristic scheduling policies.

\begin{figure}[!htbp]
	\centering
	\includegraphics[scale = 0.25]{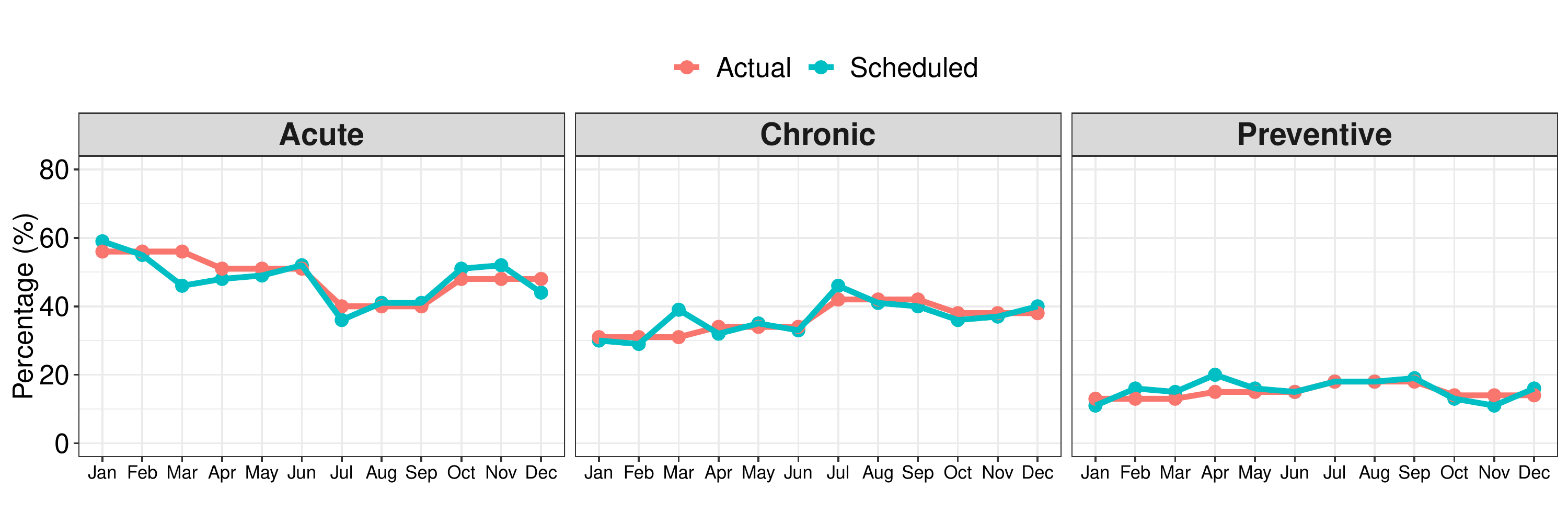}
	\caption{Actual and scheduled patterns of patient demand mix across different months at a VA primary care clinic.}
	\label{Patientmixmonthly}
\end{figure}

To consider variability in the patient mix over time, we assume that the demand mix changes across different seasons, with fewer acute patients during the summer and more during the winter, but more chronic patients during the summer and fewer during the winter. This is important because the patient types have different service times with the nurse and provider, and so any change in the demand mix can make it difficult to manage patient flow. Figure \ref{Patientmixmonthly} shows that the two-stage stochastic programming approach adjusts to changes in the patient mix and follows the pattern in different months across the planning horizon. This helps the clinic better allocate its resources and best respond to demand patterns for each patient type in each month.

\vspace{-12mm}

\subsection{Patient no-show behavior}

\begin{figure}[!htbp]
	\centering
	\includegraphics[scale = 0.3]{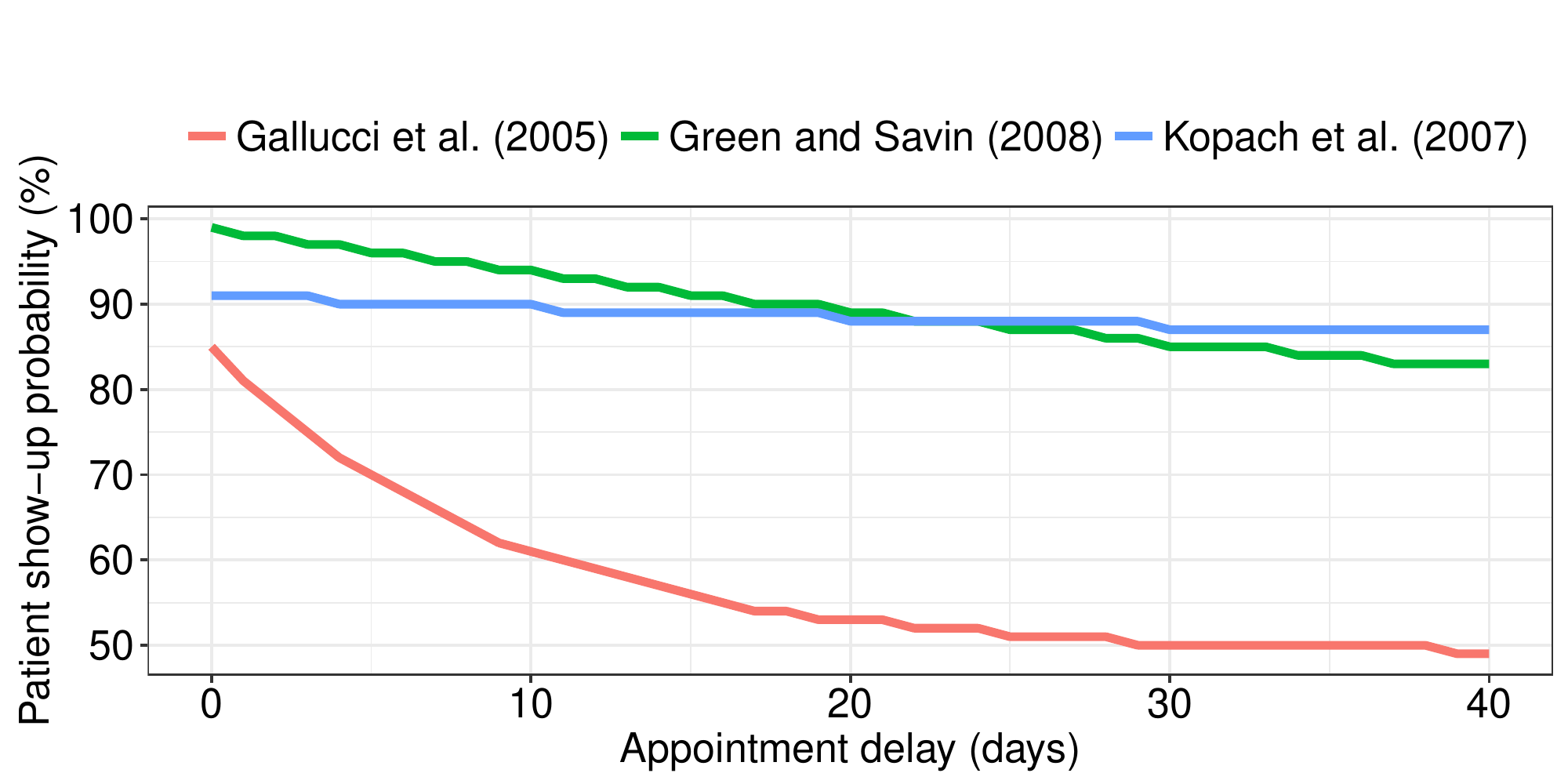}
	\caption{Patient show-up probabilities}
	\label{Patient show up probabilities}
\end{figure}

Evidence suggests that there is a higher chance that a patient will not show up when the appointment delay becomes longer \cite{festinger2002telephone} \cite{gallucci2005brief}. Three different functions have been proposed by Kopach et al. \cite{kopach2007effects}, Galluci et al. \cite{gallucci2005brief}, and Green and Savin \cite{green2008reducing} to show the relationship between appointment delays and patient show-up probabilities:

\[
p_{j} = \begin{cases}
1-p \ast (1-0.5 \ast e^{-0.017j}) & \text{Kopach et al. } \\
1-(0.51-0.36 \ast e^{-j/9}) & \text {Galluci et al. } \\
1-(0.31-0.3 \ast e^{-j/50}) & \text {Green and Savin },  
\end{cases}
\]

\noindent where the index $j$ represents the appointment delay and $p$ is the estimated patient no-show probability (in the function proposed by Kopach et al. \cite{kopach2007effects}). We assume that $p$ is equal to the average no-show probability in our study. Figure \ref{Patient show up probabilities} shows the sensitivity of patients to appointment delays under these functions. The plot shows that while the patient show-up probability converges after about 20 days under the no-show function proposed by Gallucci et al. \cite{gallucci2005brief}, the show-up probability continues to decrease in the no-show function proposed by Green and Savin \cite{green2008reducing}.

\begin{figure}[!htbp]
	\centering
	\includegraphics[scale = 0.3]{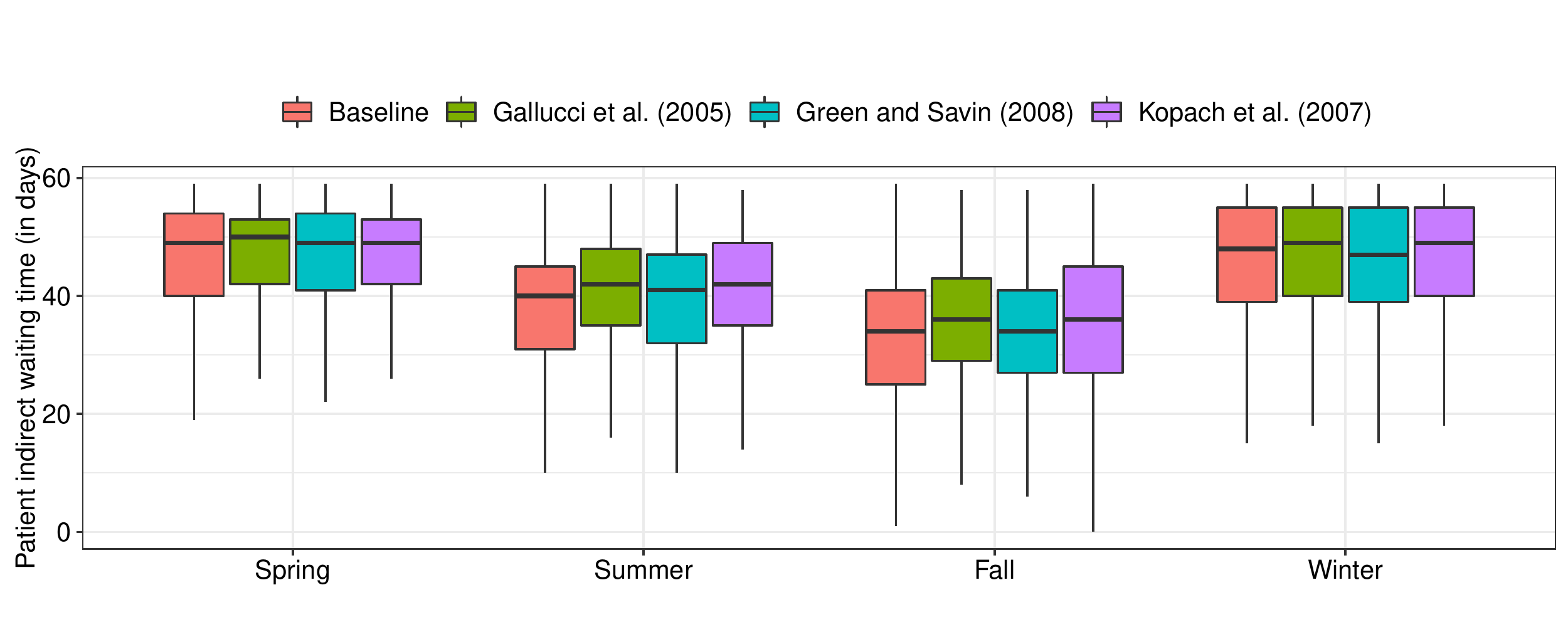}
	\caption{Patients' indirect waiting time distribution under different patient no-show behaviors across different seasons.}
	\label{IndirectWaitingDirectDailyNoshow}
\end{figure}

There is no significant difference between the various cases of constant and delay-dependent no-show functions in terms of patients' indirect waiting time. However, patients' indirect waiting time fluctuates from season to season, and the no-show behavior of patients, constant vs. delay-dependent, has a clearer impact on the patients' indirect waiting time distribution during the Summer and Fall, as shown in Figure \ref{IndirectWaitingDirectDailyNoshow}. Also, while the differences between the no-show behaviors are not significant between seasons for patients who are given same-day appointments, more patients are not given any appointment during the Spring and Winter across all no-show behaviors (around 7\% in the Spring, 0.4\% in the Summer, 1\% in the Fall, and 15\% in the Winter).

\vspace{-6mm}

\subsection{Provider capacity}

Primary care providers need to devote their time and capacity to different activities such as clinical work, teaching, research, surgery, and more. This may result in a provider having less than the expected time for clinical care and the office practice. The full-time equivalent (FTE) for an activity is the percentage of a provider's time that is spent on that activity. As shown in the base case scenario, the primary care provider is assigned one appointment scheduling session per week, in the morning or in the afternoon, for non-clinical activities. In this section, we analyze the impact of the provider's availability on patients' waiting times.

\begin{figure}[!htbp]
	\centering
	\includegraphics[scale = 0.3]{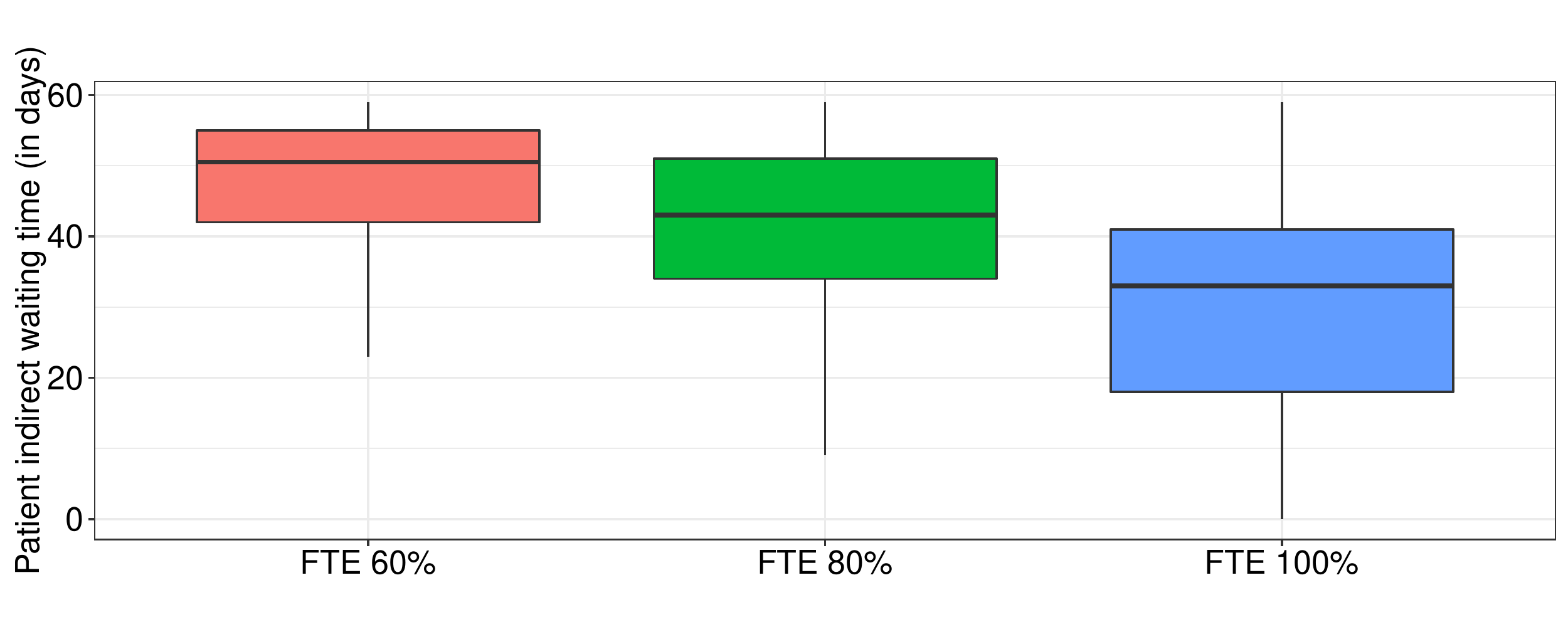}
	\caption{Patients' indirect waiting time distribution under different provider FTEs.}
	\label{IndirectWaitingDirectProviderType}
\end{figure}

Figure \ref{IndirectWaitingDirectProviderType} represents the distribution of patients' indirect waiting times for different clinical FTEs. Assigning more time to clinical activities will help more patients get appointments with the provider. Our analysis shows that 74\% of patients can get an appointment with the provider if 60\% of the provider's time is assigned to clinical activities, while 3\% of the patients can see the provider on the same day, and 100\% of the patients can be given an appointment if all of the provider's time is devoted to visiting patients. 

\subsection{Nurse and provider service time distributions}

The nurses' and providers' expertise impact the flow of patients, so service time can vary from nurse to nurse and from provider to provider. The provider might also be a provider assistant (PA) or nurse practitioner (NP) instead of an actual provider, in which case the patient panel characteristics could be different. For example, in the Veterans Health Administration, NPs visit more women (10\% of patients) than do PAs (6.7\%) and providers (6.6\%), and providers and NPs have more patients from minority groups (21\% and 20\%, respectively) in their panel, compared to PAs (18\%) \cite{morgan2012characteristics}. In addition, medical complexity and the number of new patients can cause more variability in nurse and provider service time.  

The effect of variability in the CV of the nurse and provider service time is evaluated by considering three cases: change only in the nurse's CV, change only in the provider's CV, and change in both the nurse's and the provider's CV. In the test cases, the standard deviation is multiplied by either 0.8 or 1.2 to give lower or higher CVs.

\begin{figure}[!htbp]
	\centering
	\includegraphics[scale = 0.3]{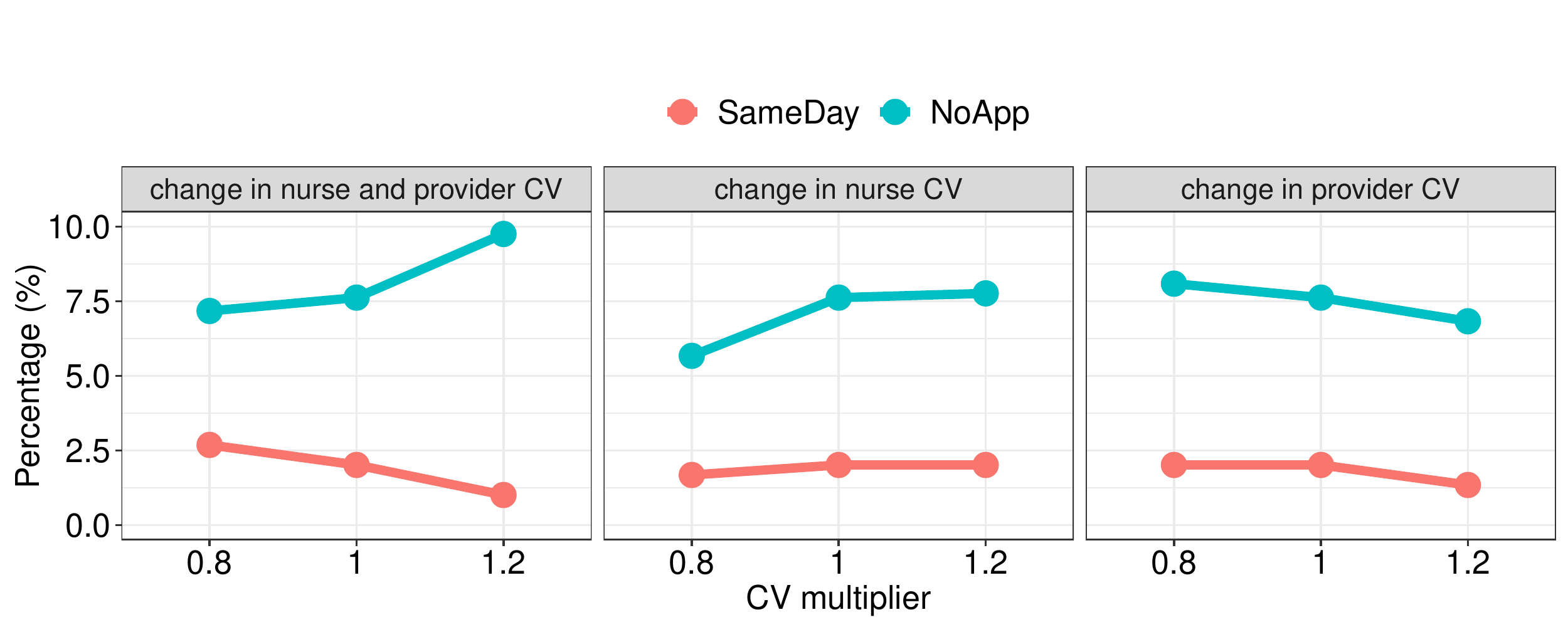}
	\caption{Percentage of patients who are given same-day appointments (red) and  percentage of patients who are not given any appointment (blue) under different values of the CV multiplier for the nurse and the provider.}
	\label{IndirectWaitingDirectServiceTimeSameDay}
\end{figure}

\begin{figure}[!htbp]
	\centering
	\includegraphics[scale = 0.3]{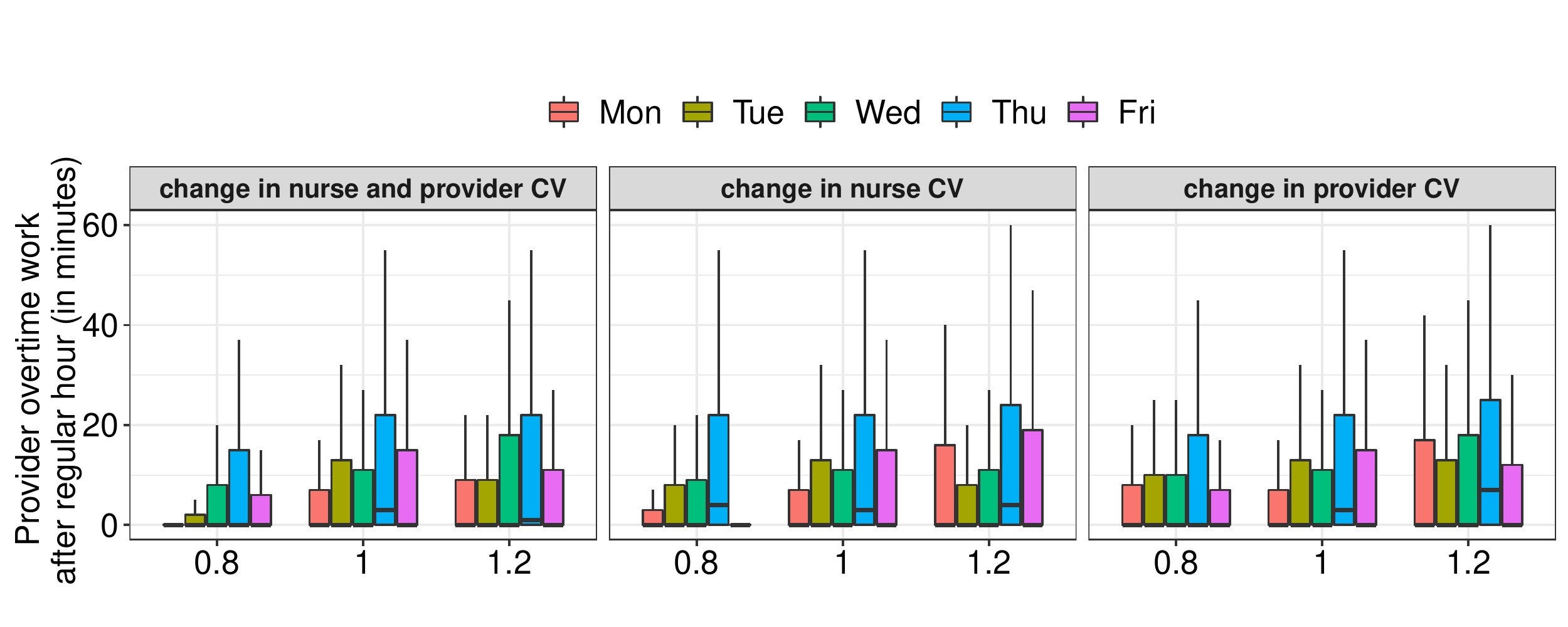}
	\caption{Distribution of provider overtime work after regular hours under different values of the CV multiplier for the nurse and the provider each day of the week.}
	\label{ProviderOvertimeServiceTime}
\end{figure}

Although the patients' indirect waiting time distribution in each season and the percentage of patients who receive same-day appointments are minimally different for the three cases, the percentage of patients who do not receive any appointment increases as both the nurse and provider have higher CV values (Figure \ref{IndirectWaitingDirectServiceTimeSameDay}). Moreover, Figure \ref{ProviderOvertimeServiceTime} shows that a change in the provider's CV has the highest impact on the amount of time that the provider has to work after regular hours to serve patients. 

\vspace{-12mm}

\subsection{Appointment slot perishability}

\begin{figure}[!htbp]
	\centering
	\includegraphics[scale = 0.3]{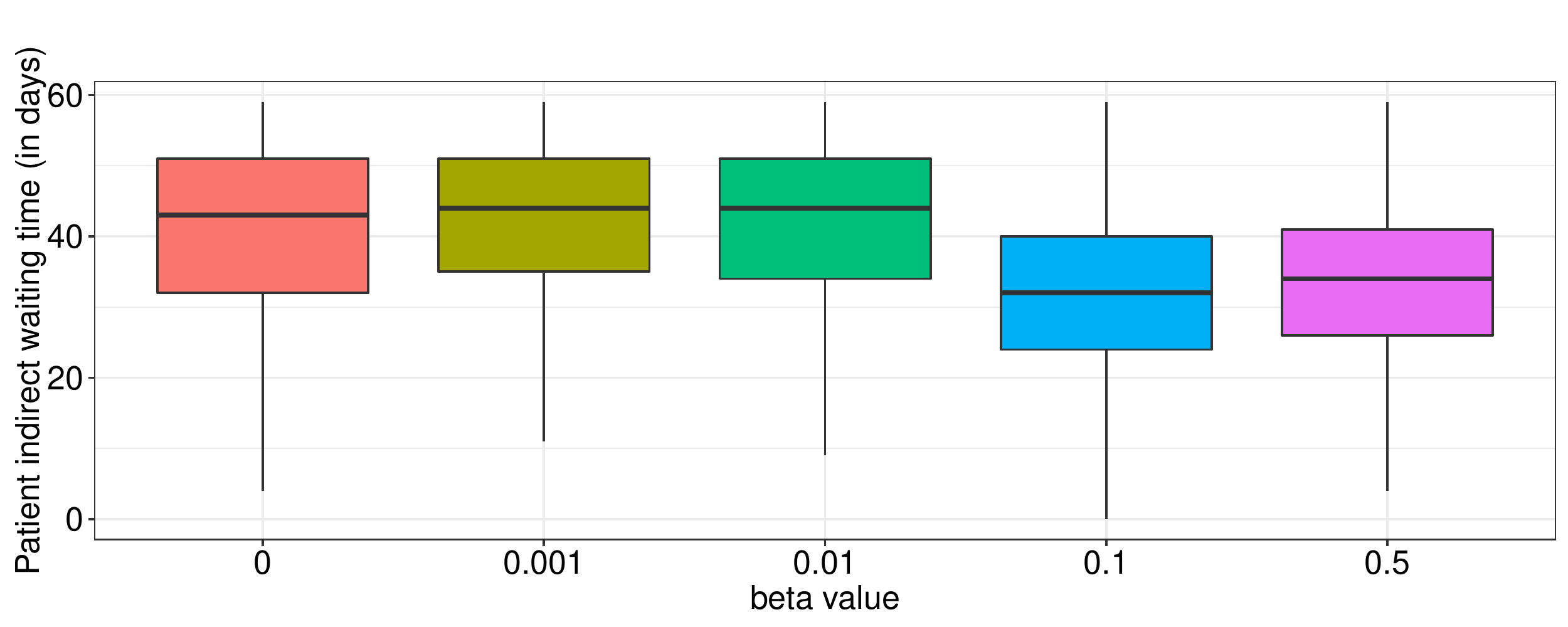}
	\caption{Patients' indirect waiting time distribution under different values of beta in the index policy.}
	\label{IndirectWaitingIndexPolicy}
\end{figure}

In the computational experiments, we assumed that the call center uses the index policy to simulate the appointment scheduling in practice. One major drawback of the heuristic index policy is that it does not consider the ``perishability'' of open slots, i.e., the index policy ranks the days with open appointments without considering how far they are from the patient's desired date. To overcome this drawback, we propose the following modified index for each appointment slot in the scheduling horizon:

\begin{equation}
I(j) = c_j*e^{(\beta*(DD - date_{j})},
\end{equation}

\noindent where $c_{j}$ is the remaining capacity in appointment slot $j$, $DD$ stands for a patient's desired date, and $date_{j}$ is the date of the appointment slot $j$. As $\beta$ decreases, this converges to the simple index policy. Figure \ref{IndirectWaitingIndexPolicy} shows that when this approach is used, there is significant improvement in patients' indirect waiting time as the value of $\beta$ increases. This means that if the scheduler considers how close an open appointment is to the patient's desired date instead of the capacity of the open appointment, the patient can visit the primary care provider sooner. In addition, a smaller percentage of patients will fail to be given an appointment with the provider as the value of $\beta$ increases (from 7.62\% when $\beta = 0$ to 3.72\% when $\beta = 0.5$), while the percentage of patients who get same-day appointments will not change significantly (from 2\% when $\beta = 0$ to 3\% when $\beta = 0.5$).  

\vspace{-9mm}

\section{Conclusions and future work} \label{Conclusion}

Although primary care is considered to be patients' first point of contact with a healthcare system, patients often suffer from significant delays in obtaining appointments. A well-established appointment scheduling system can help clinics reduce patients' indirect waiting time while also improving patient flow within the clinic. Clinic managers have to handle multiple issues when scheduling patient appointments. While different patient types have different complexities, there are uncertainties in the pattern of calls for appointments and patients' willingness to wait for an appointment. Some patients may call in advance to book appointments, while others ask for same-day appointments. This study proposes solving the appointment scheduling problem with a two-stage stochastic programming model integrated with a simulation model to minimize patients' indirect waiting time for an appointment while maintaining a patient flow in the clinic that is within the constraints of acceptable patient and practitioner performance thresholds. The model proposes a patient scheduling template for the call center in order to help the clinic manager reduce patient appointment delays and scheduling errors and improve the efficiency of resource allocation.

We use a numerical case study inspired by a real-world healthcare system to validate our proposed approach over heuristic approaches suggested in the literature. Our model also performs better when the provider is assigned greater capacity for other clinical activities as well as when appointment slot perishability is considered in scheduling the patient appointments. Determining the optimal schedule for re-running the optimization model based on changes in the uncertainties is an important avenue for future research. Moreover, while our two-stage stochastic programming model determines the number of each patient type that can be scheduled in each appointment slot, it does not provide any guidance regarding when the call center should offer each open appointment to a patient. Integrating the call center scheduling process with our two-stage stochastic programming approach could result in further improving patients' indirect waiting times.  

\begin{acknowledgements}
	We thank the U.S. Veteran's Health Administration for sponsoring part of this research.
\end{acknowledgements}

%
\section*{Conflict of interest}
The authors declare that they have no conflict of interest.

\bibliographystyle{spmpsci}      
\bibliography{Bio}   

%
%

\end{document}